\documentclass{amsart}

\usepackage{amsmath}
\usepackage{amscd}
\usepackage{amssymb}

\newcommand{\bP}{{\mathbb P}}

\newcommand{\bZ}{{\mathbb Z}}

\newcommand{\cC}{{\mathcal C}}
\newcommand{\cE}{{\mathcal E}}
\newcommand{\cM}{{\mathcal M}}
\newcommand{\cO}{{\mathcal O}}
\newcommand{\cP}{{\mathcal P}}

\newcommand{\cW}{{\mathcal W}}
\newcommand{\Mbar}{\overline{\cM}}

\DeclareMathOperator{\Aut}{Aut}

\newtheorem{conjecture}{Conjecture}

\newtheorem{theorem}{Theorem}[section]
\newtheorem{theorem/definition}{Theorem/Definition}[section]

\newtheorem{proposition}{Proposition}[section]
\newtheorem{lemma}{Lemma}[section]

\newtheorem{Conjecture}{Conjecture}

\theoremstyle{remark}
\newtheorem{remark}{Remark}[section]

\theoremstyle{definition}
 \newtheorem{example}{Example}[section]

\begin{document}

\title
{A Conjecture ON HODGE INTEGRALS}
\author{Jian Zhou}
\address{Department of Mathematical Sciences\\Tsinghua University\\Beijing, 100084, China}
\email{jzhou@math.tsinghua.edu.cn} \footnote{This research
is partially supported by research grants from NSFC and Tsinghua
University.}
\begin{abstract}
We propose a conjectural formula expressing the generating series of some Hodge integrals
in terms of representation theory of Kac-Moody algebras.
Such generating series appear in calculations of Gromov-Witten invariants by localization techniques.
It generalizes a formula conjectured by Mari\~no and Vafa,
recently proved in joint work with Chiu-Chu Melissa Liu and Kefeng Liu.
Some examples are presented.
\end{abstract}
\maketitle

\section{Introduction}

An integral of the form
\begin{eqnarray} \label{eqn:Hodge}
\int_{\Mbar_{g, n}} \psi_1^{i_1} \cdots\psi_n^{i_n}
\lambda_1^{j_1} \cdots \lambda_g^{j_g}
\end{eqnarray}
is called a Hodge integral.
Here $\psi_i$ and $\lambda_j$ are Chern classes of some naturally defined
vector bundles on the moduli space $\Mbar_{g, n}$.
One way to evaluate such integrals is to reduce them to the integrals of the form
\begin{eqnarray} \label{eqn:Psi}
\int_{\Mbar_{g, n}} \psi_1^{i_1} \cdots\psi_n^{i_n},
\end{eqnarray}
which are covered by the Witten conjecture/Kontsevich theorem \cite{Wit1, Kon1},
and an algorithm is available to do it automatically \cite{Fab}.
This has the disadvantage of having to be done genus by genus.

Localization methods are the most powerful mathematical techniques in the computations of
Gromov-Witten invariants (cf. \cite{Kon, Gra-Pan}).
Hodge integrals naturally appear in such computations
in positive genus.
As suggested by physicists,
it is crucial to take generating series in all genera
and get a closed expression
from which one can extract important integral invariants \cite{Gop-Vaf1, Gop-Vaf3}.
For example,
the following generating series of Hodge integrals appear in the localization
calculations for toric Fano surfaces \cite{Zho4}:
\begin{equation*} \label{eqn:CZ}
\begin{split}
& G_{\mu^+, \mu^-}(\lambda; x, y) \\
& =  - \frac{\sqrt{-1}^{l(\mu^+)+l(\mu^-)}}{z_{\mu^+} \cdot z_{\mu^-}} \\
&\cdot  \sum_{g \geq 0} \lambda^{2g-2} \int_{\Mbar_{g, l(\mu^+)+l(\mu^-)}}
\frac{\Lambda_{g}^{\vee}(x)\Lambda^{\vee}_{g}(y)\Lambda_{g}^{\vee}(-x - y)}
{\prod_{i=1}^{l(\mu^+)} \frac{x}{\mu_i^+} \left(\frac{x}{\mu^+_i} - \psi_i\right)
\prod_{j=1}^{l(\mu^-)} \frac{y}{\mu_i^-}\left(\frac{y}{\mu^-_j} - \psi_{j+l(\mu^+)}\right)} \\
& \cdot \left[xy(x+y)\right]^{l(\mu^+)+l(\mu^-)-1}
\cdot \prod_{i=1}^{l(\mu^+)} \frac{\prod_{a=1}^{\mu^+_i-1}
\left( \mu^+_iy + a x\right)}{\mu_i^+! x^{\mu_i^+-1}}
\cdot \prod_{i=1}^{l(\mu^-)} \frac{\prod_{a=1}^{\mu^-_i-1}
\left( \mu_i^- x + a y\right)}{\mu_i^-! y^{\mu_i^--1}}.
\end{split} \end{equation*}
Here $\mu^+$ and $\mu^-$ are two partitions,
one of which might be empty,
and
$$\Lambda_g^{\vee}(x) = \sum_{i=0}^g (-1)^i\lambda_ix^{g-i}.$$
When $\mu^-$ is empty,
the generating series becomes:
\begin{equation} \label{eqn:CMV}
\begin{split}
\cC_{\mu}(\lambda; x, y) &=  \sum_{g \geq 0} \lambda^{2g-2} \int_{\Mbar_{g, l(\mu^+)}}
\frac{\Lambda_{g}^{\vee}(x)\Lambda^{\vee}_{g}(y)\Lambda_{g}^{\vee}(-x - y)}
{\prod_{i=1}^{l(\mu^+)} \frac{x}{\mu_i^+} \left(\frac{x}{\mu^+_i} - \psi_i\right)} \\
& \cdot \left[xy(x+y)\right]^{l(\mu^+)-1}
\cdot \prod_{i=1}^{l(\mu^+)} \frac{\prod_{a=1}^{\mu^+_i-1}
\left( \mu^+_iy + a x\right)}{\mu_i^+! x^{\mu_i^+-1}},
\end{split} \end{equation}
which appear in the localization calculations
for $\cO(-1) \oplus \cO(-1) \to \bP^1$,
both in the open string context \cite{Kat-Liu} and the closed string context \cite{Zho2,Zho5,Zho6}).

String theory provides an unexpected link between Hodge integrals
with other branches of mathematics.
Mari\~no and Vafa \cite{Mar-Vaf} conjectured by string duality a formula relating (\ref{eqn:CMV})
to Wess-Zumino-Witten (WZW) theory.
The explicit form of the Mari\~no-Vafa conjecture appears to involve only the representation
theory of symmetric groups,
but the physical procedure of deriving it involves many of the deep ideas developed
by string theorists,
e.g. the relationship between Chern-Simons theory with link variants and representation
theory of Kac-Moody algebras.

Some special cases of this formula were first proved in \cite{Zho1,Zho2}.
In \cite{Zho3} the author proposed the following approach to prove the Mari\~no-Vafa formula:
One first shows both sides of the formula satisfy the cut-and-join equation,
then shows that they have the same initial values.
The relevant combinatorial issues are also dealt with in \cite{Zho3},
including the establishment of the cut-and-join equation for the combinatorial side
of the Mari\~no-Vafa formula and the identification of the initial values.
In joint work with Chiu-Chu Melissa Liu and Kefeng Liu \cite{LLZ1, LLZ2},
we also also establish the cut-and-join equation for the geometric side of the Mari\~no-Vafa formula
hence completing the proof of this formula. (See \cite{Oko-Pan}) for a differen approach.)
Some applications of this formula to Hodge integrals can be found in \cite{LLZ3}.

In \cite{Zho4} we show how to apply the Marino-Vafa
formula to calculate the BPS numbers in local $\bP^2$ and $\bP^1 \times \bP^1$ geometry.
In \cite{Zho5} we announce a proof of a conjecture by Iqbal \cite{Iqb} relating the local
geometry of $\cO(-1) \oplus \cO(-1) \to \bP^1$ to WZW theory.
The efforts to generalize these works to more general local toric surface geometries
lead us to this work.
In this work we make the following:

\begin{Conjecture}
\begin{equation} \label{eqn:Conj}
\begin{split}
& \exp \left(
\sum_{(\mu^+, \mu^-) \in \cP^2} G_{\mu^+, \mu^-}(x, y)p^+_{\mu^+}p^-_{\mu^-}\right) \\
= & \sum_{|\mu^{\pm}|=|\nu^{\pm}|}
\frac{\chi_{\nu^+}(\mu^+)}{z_{\mu^+}} \frac{\chi_{\nu^-}(\mu^-)}{z_{\mu^-}}
e^{\sqrt{-1}(\kappa_{\nu^+} \frac{y}{x}  + \kappa_{\nu^-} \frac{x}{y})\lambda/2}
\cW_{\nu^+, \nu^-} p^+_{\mu^+}p^-_{\mu^-}.
\end{split} \end{equation}
\end{Conjecture}

See Section \ref{sec:Pre} and Sections \ref{sec:Conjecture} for notations.

We will propose a proof of this conjecture along the same lines of
our proof of the Mari\~no-Vafa formula.
I.e.,
I will show that
a certain cut-and-join equation is satisfied by the right-hand side of (\ref{eqn:Conj}),
and both sides have the same initial values.
In a forthcoming paper with Chiu-Chu Melissa Liu and Kefeng Liu \cite{LLZ4},
we will use localization method to prove the left-hand side
of (\ref{eqn:Conj}) also satisfies the same equation,
hence complete the proof of (\ref{eqn:Conj}).
The applications including the proof of Iqbal's conjecture for general local toric geometry
will be presented in \cite{Zho6}.

The rest of this paper is arranged as follows.
In Section \ref{sec:Pre} we recall some preliminary facts from
Wess-Zumino-Witten theory,
which corresponds to the representation theory
of integrable highest weight representations
of Kac-Mody algebras in mathematical literature.
We recall the Mari\~no-Vafa formula and present our conjecture in Section \ref{sec:Conjecture}.
A proof based on the cut-and-join equation is proposed in this section.
We recall some facts about skew Schur functions and prove an orthogonality relation for them
in Section \ref{sec:Skew}.
A reformulation of $\cW_{\nu^+, \nu^-}$ in terms of skew Schur functions
is given in Section \ref{sec:RefW}.

We study the initial value problem for the relevant cut-and-join equation in Section \ref{sec:Initial}.
Some examples are presented in Sections \ref{sec:Examples}.

\section{Preliminaries}
\label{sec:Pre}

\subsection{Partitions}
We use Macdonald's book \cite{Mac} as our reference.
A partition of a positive integer $d$ is a sequence of integers $n_1 \geq n_2 \geq \dots \geq n_l > 0$
such that $n_1 + \dots + n_l = d$.
We write
\begin{align*}
|\eta| & = d, &
l(\eta) & = l.
\end{align*}
Denote by $m_j(\eta)$ the number of $j$'s among $n_1, \dots, n_l$.
Each partition $\eta$ of $d$ corresponds to a conjugacy class $C_{\eta}$ of $S_d$.
Denoted by $C_{(2)}$ the conjugacy class of transpositions.
The number of elements in $C_{\eta}$ is
$$|C_{\eta}| = \frac{d!}{z_{\eta}},$$
where
$$z_{\eta} = \prod_j m_j(\eta)!j^{m_j(\eta)} .$$
Denote by $(-1)^g$ the sign of an element in $S_d$.
It is easy to see that
$$(-1)^g = (-1)^{|\eta| - l(\eta)},$$
for $g \in C_{\eta}$.

Another way of representing a partition is by its Young diagram.
The Young diagram of $\eta$ has $m_j(\eta)$ rows of boxes of length $j$.
The partition corresponding to the transpose of the Young diagram of $\eta$
will be denoted by $\eta'$.
The number of squares in the $i$-th row of $\eta'$ will be written as $\eta_i'$.
For any square $e \in \eta$,
denote by $h(e)$ its hook length.

Each partition $\lambda$ corresponds to an irreducible representation $R_{\lambda}$ of $S_d$.
For example,
$\chi_{(d)}$ corresponds to the trivial representation,
$\chi_{(1^d)}$ corresponds to the sign representation.
The value of the character $\chi_{R_{\lambda}}$ on the conjugacy class $C_{\eta}$ is
denoted by $\chi_{\lambda}(\eta)$.

Sometimes we will also need the partition of $0$.
Denote $\cP$ the set of all partitions of nonegative integers,
by $\cP_+$ the set of all partitions of positive integers.
Denote by $\cP^2 = \cP \times \cP$ the set of pairs of partitions,
and by $\cP^2_+$ the set of pairs $(\mu^+, \mu^-) \in \cP \times \cP$
such that $\mu^+$ and $\mu^-$ are not both partitions of $0$.

\subsection{Symmetric functions}
Denote by $p_n$ the $n$-th symmetric power function of degree $n$.
For a partition $\mu$,
define
$$p_{\mu} = \prod_i p_i^{m_i(\mu)}.$$
Similarly define $e_{\mu}$,
where $e_n$ is the $n$-th elementary symmetric function.
Then $\{p_{\mu}\}_{\mu \in \cP}$ and $\{e_{\mu}\}_{\mu \in \cP}$
are bases of the space of symmetric functions.
Another basis is given by the Schur functions $\{s_{\mu}\}_{\mu \in \cP}$.
They are related to $e_n$'s by the famous Jacobi-Trudy formula:
$$s_{\mu}(\cE(t)) = \det (e_{\mu'_i-i +j}),$$
where $\cE(t) = \sum_{i \geq 0} e_i t^i$.

\subsection{Modular property of integrable highest weight representations}
For a fixed integer $k$,
there are only finitely many integrable highest weight representations of level $k$
of an affine Kac-Moody algebra up to equivalence.
Denote their characters by $\chi_0(\tau), \dots, \chi_n(\tau)$.
Then there are holomorphic functions $S_{ij}(\tau)$,
such that
\begin{eqnarray*}
&& \chi_i(- \frac{1}{\tau}) = \sum_j S_{ij}(\tau) \chi_j(\tau).
\end{eqnarray*}
From this one can construction a representation of a double covering of
$SL(2, \bZ)$ (cf. \cite{Ver}).

The $U(N)$ WZW theory are associated with integrable highest weight representations
of level $k+N$ \cite{Ver, Wit}.
Such representations are indexed by partitions.
The matrix elements of $S^{-1}$ is denoted by $W_{\mu, \nu}$,
where $\mu, \nu$ are partitions which might be empty.
They are given by the Morton-Lukac formula \cite{Mor-Luk, Mar-Vaf, Aga-Mar-Vaf, Iqb}.
For us, as in the above references,
certain leading term in the large $N$ expansion denoted by $\cW_{\mu, \nu}$ will be used.
They can be computed as follows:
\begin{eqnarray}
&& \cW_{\mu, \nu} = q^{l(\nu)/2} \cW_{\mu} \cdot s_{\nu}(\cE_{\mu}(t)),
\end{eqnarray}
where
\begin{eqnarray}
&& \cW_{\mu} = q^{\kappa_{\mu}/4}\prod_{1 \leq i < j \leq l(\mu)}
\frac{[\mu_i - \mu_j + j - i]}{[j-i]}
\prod_{i=1}^{l(\mu)} \prod_{v=1}^{\mu_i} \frac{1}{[v-i+l(\mu)]}, \\
&& \cE_{\mu}(t) = \prod_{j=1}^{l(\mu)} \frac{1+q^{\mu_j-j}t}{1+q^{-j}t}
\cdot \left(1 + \sum_{n=1}^{\infty}
\frac{t^n}{\prod_{i=1}^n (q^i-1)}\right).
\end{eqnarray}
As usual,
$$[m] = q^{m/2} - q^{-m/2}.$$

It is not obvious from the above expression that we actually have:
\begin{eqnarray} \label{eqn:Symetric}
&& \cW_{\mu, \nu} = \cW_{\nu, \mu}.
\end{eqnarray}
In \S \ref{sec:RefW} we will reformulate $\cW_{\mu, \nu}$ to make this symmetry manifest.

\section{The Conjecture}

\label{sec:Conjecture}

\subsection{Mari\~no-Vafa formula}
To motivate our conjecture,
we first recall this remarkable formula.
Consider the following series:
\begin{eqnarray*}
\cC_{\eta}(\tau, \lambda) & = &
- \frac{\sqrt{-1}^{l(\eta)}}{z_{\eta}} \cdot
(\tau(\tau+1))^{l(\eta)-1}  \cdot \prod_{i=1}^{l(\eta)} \frac{\prod_{j=1}^{\eta_i-1}(j+\eta_i\tau)}{\eta_i!}\\
&& \cdot\sum_{g \geq 0} \lambda^{2g-2+l(\eta)}
\int_{\Mbar_{g, l(\eta)}} \frac{\Lambda^{\vee}_g(1)\Lambda^{\vee}_g(-\tau-1)\Lambda^{\vee}_g(\tau)}
{\prod_{i=1}^{l(\eta)} \frac{1}{\eta_i}(\frac{1}{\eta_i} -\psi_i)}, \\
\cC(\tau, \lambda, p) & = & \sum_{\eta} \cC_{\eta}(\tau, \lambda)p_{\eta}, \\
\cC(\tau, \lambda, p)^{\bullet} & = & \exp \cC(\tau, \lambda, p).
\end{eqnarray*}
Then the Mari\~no-Vafa formula can be written as follows:
\begin{equation} \label{eqn:Mar-Vaf2}
\begin{split}
& \cC_{\eta}(\tau,\lambda)
= \sum_{n \geq 1}
\frac{(-1)^{n-1}}{n}\sum_{\cup_{i=1}^n \eta_i = \eta}
\prod_{i=1}^n \sum_{|\rho_i|=|\eta_i|}
\frac{\chi_{\rho_i}(\eta_i)}{z_{\eta_i}}
 \cdot  e^{\sqrt{-1}\tau\kappa_{\rho_i}\lambda/2} \cdot \cW_{\rho_i},
\end{split} \end{equation}
or equivalently,
$$\cC(\tau, \lambda, p)^{\bullet}
= \sum_{|\rho|=|\eta| \geq 0}
\frac{\chi_{\rho}(\eta)}{z_{\eta}}
 \cdot  e^{\sqrt{-1}\tau\kappa_{\rho}\lambda/2} \cdot \cW_{\rho}p_{\eta}.$$
This was first explicitly written down in \cite{Zho3} and first proved in \cite{LLZ1, LLZ2}.

\subsection{Our conjecture}
Consider the following generating series of Hodge integrals:
\begin{eqnarray*}
&& G_{\mu^+, \mu^-}(x, y) \\
& = & - \frac{\sqrt{-1}^{l(\mu^+)+l(\mu^-)}}{z_{\mu^+} \cdot z_{\mu^-}} \\
&& \cdot \sum_{g \geq 0} \lambda^{2g-2} \int_{\Mbar_{g, l(\mu^+)+l(\mu^-)}}
\frac{\Lambda_{g}^{\vee}(x)\Lambda^{\vee}_{g}(y)\Lambda_{g}^{\vee}(-x - y)}
{\prod_{i=1}^{l(\mu^+)} \frac{x}{\mu_i^+} \left(\frac{x}{\mu^+_i} - \psi_i\right)
\prod_{j=1}^{l(\mu^-)} \frac{y}{\mu_i^-}\left(\frac{y}{\mu^-_j} - \psi_{j+l(\mu^+)}\right)} \\
&& \cdot \left[xy(x+y)\right]^{l(\mu^+)+l(\mu^-)-1}
\cdot \prod_{i=1}^{l(\mu^+)} \frac{\prod_{a=1}^{\mu^+_i-1}
\left( \mu^+_iy + a x\right)}{\mu_i^+! x^{\mu_i^+-1}}
\cdot \prod_{i=1}^{l(\mu^-)} \frac{\prod_{a=1}^{\mu^-_i-1}
\left( \mu_i^- x + a y\right)}{\mu_i^-! y^{\mu_i^--1}}.
\end{eqnarray*}

Recall the following easy result from \cite{Tia-Zho}:

\begin{lemma} \label{lm:Degree}
We have the following identity:
\begin{eqnarray*}
&& \int_{\Mbar{g, n}}
\frac{\prod_{j=1}^r \Lambda^{\vee}_g(a_j t)}
{\prod_{k=1}^n\frac{t}{d_k}\left(\frac{t}{d_k}-\psi_k\right)}
 =  t^{(r-3) g - 3n  + 3}
\int_{\Mbar_{g, n}} \frac{\prod_{j=1}^r \Lambda^{\vee}_g(a_j)}
{\prod_{k=1}^n \frac{1}{d_k}\left(\frac{1}{d_k}-\psi_k\right)}.
\end{eqnarray*}
\end{lemma}

By Lemma \ref{lm:Degree},
\begin{eqnarray*}
&& G_{\mu^+, \mu^-}(x, y) \\
& = &  -\frac{\sqrt{-1}^{l(\mu^+)+l(\mu^-)}}{|\Aut(\mu^+)| \cdot |\Aut(\mu^-)|}
\lambda^{l(\mu^+)+l(\mu^-)-2} \\
&& \cdot \sum_{g \geq 0} \lambda^{2g} \int_{\Mbar_{g, l(\mu^+)+l(\mu^-)}}
\frac{\Lambda_{g}^{\vee}(1)\Lambda^{\vee}_{g}(\tau)\Lambda_{g}^{\vee}(-1 - \tau)}
{\prod_{i=1}^{l(\mu^+)} \frac{1}{\mu_i^+} \left(\frac{1}{\mu^+_i} - \psi_i\right)
\prod_{j=1}^{l(\mu^-)} \frac{\tau}{\mu^-_j} \left(\frac{\tau}{\mu^-_j} - \psi_{j+l(\mu^+)}\right)} \\
&& \cdot \left[\tau(1+\tau)\right]^{l(\mu^+)+l(\mu^-)-1}
\cdot \prod_{i=1}^{l(\mu^+)} \frac{\prod_{a=1}^{\mu^+_i-1}
\left( \mu_i^+ \tau + a \right)}{\mu_i^+!}
\cdot \prod_{i=1}^{l(\mu^-)} \frac{\prod_{a=1}^{\mu^-_i-1}
\left(\frac{\mu_i^-}{\tau} + a \right)}{\mu_i^-!}.
\end{eqnarray*}
Here
$$\mu= \frac{y}{x}.$$
Hence we will write $G_{\mu^+, \mu^-}(x, y)$ as $G_{\mu^+, \mu^-}(\tau)$.
Since
$$G_{\mu^+, \mu^-}(x, y) = G_{\mu^-, \mu^+}(y, x),$$
we clearly have
\begin{eqnarray} \label{eqn:Inverse}
&& G_{\mu^+, \mu^-}(\tau) = G_{\mu^-, \mu^+}(\frac{1}{\tau}).
\end{eqnarray}
Hence
\begin{eqnarray*}
&& G_{\mu^+, \mu^-}(-\frac{1}{\tau}) = G_{\mu^-, \mu^+}(-\tau).
\end{eqnarray*}

Now we can state the main subject of this paper.

\begin{conjecture}
We have the following identity
\begin{equation} \label{eqn:Conjecture}
\begin{split}
& \exp \left(
\sum_{(\mu^+, \mu^-) \in \cP_+^2} G_{\mu^+, \mu^-}(\tau)p^+_{\mu^+}p^-_{\mu^-}\right) \\
= & \sum_{|\mu^{\pm}|=|\nu^{\pm}| \geq 0}
\frac{\chi_{\nu^+}(\mu^+)}{z_{\mu^+}} \frac{\chi_{\nu^-}(\mu^-)}{z_{\mu^-}}
e^{\sqrt{-1}(\kappa_{\nu^+} \tau  + \frac{\kappa_{\nu^-}}{\tau})\lambda/2}
\cW_{\nu^+, \nu^-} p^+_{\mu^+}p^-_{\mu^-}.
\end{split} \end{equation}
\end{conjecture}

\subsection{Cut-and-join equation}
Write the right-hand side of (\ref{eqn:Conjecture}) as $R(\tau)^{\bullet}$.
Here we suppress $\lambda, p^{\pm}$ from the notation for simplicity.
Then by the method of \cite{Zho3} it is straightforward to see that
assuming (\ref{eqn:Conjecture}) one has

\begin{theorem}
The following equation is satisfied:
\begin{equation} \label{eqn:CutJoin}
\begin{split}
\frac{\partial}{\partial \tau} R(\tau)^{\bullet} = \frac{1}{2}(C^++J^+)R(\tau)^{\bullet}
- \frac{1}{2\tau^2}(C^-+J^-)R(\tau)^{\bullet},
\end{split}
\end{equation}
where
\begin{align*}
C^{\pm} & = \sum_{i,j} p_i^{\pm}p_j^{\pm} \frac{\partial}{\partial p^{\pm}_{i+j}}, &
J^{\pm} & = \sum_{i,j} p_{i+j}^{\pm} \frac{\partial p_i^{\pm}}{\partial p^{\pm}_j}.
\end{align*}
\end{theorem}

In Section \ref{sec:Initial} below,
we will show both sides of (\ref{eqn:Conjecture}) have the same values at $\tau = -1$.
In \cite{LLZ4} we will show the left-hand side of (\ref{eqn:Conjecture}) satisfies the same
equation as above,
hence prove (\ref{eqn:Conjecture}).

\section{Skew Schur Polynomials and Skew Schur Functions}
\label{sec:Skew}

In this section we recall some basic definitions and facts about skew Schur functions.
We will also prove an orthogonality relation that will be used later.

\subsection{Definition}
The {\em skew Schur function} $s_{\mu/\nu}$ is defined to be the symmetric function
such that
$$\langle s_{\mu/\nu}, s_{\rho}\rangle
= \langle s_{\mu}, s_{\nu}s_{\rho}\rangle.$$
Equivalently,
suppose
$$s_{\nu}s_{\rho} = \sum_\mu c^{\mu}_{\nu\rho} s_{\mu},$$
then
$$s_{\mu/\nu} = \sum_{\rho} c^{\mu}_{\nu\rho}s_{\rho}.$$
Note $s_{\mu/\nu}$ is homogeneous of degree $|\mu|-|\nu|$.

\subsection{Determinatal formula}
Recall \cite{Mac}
\begin{eqnarray} \label{eqn:SkewDet}
&& s_{\mu/\nu} = \det (h_{\lambda_i - \mu_j-i+j})_{1 \leq i, j \leq n}
= \det (e_{\lambda_i'-\mu_j'-i+j})_{1 \leq i, j \leq m},
\end{eqnarray}
where $n \geq l(\lambda)$, $m \geq l(\lambda')$.
It follows that $s_{\mu/\nu} = 0$ unless $\mu \subset \nu$,
and
$$\omega(s_{\mu/\nu}) = s_{\mu'/\nu'}.$$
Hence we have
\begin{eqnarray*}
&& c^{\mu'}_{\nu'\rho'} = c^{\mu}_{\nu\rho}, \\
&& \sum_{\rho} c^{\mu}_{\nu\rho} s_{\rho'} = s_{\mu'/\nu'}.
\end{eqnarray*}

\subsection{An orthogonality relation}
The following result seems to be new.

\begin{lemma}
We have for any two partitions $\mu$ and $\nu$
the following identity:
\begin{eqnarray} \label{eqn:Orthogonality}
&& \sum_{\rho} (-1)^{|\rho|} s_{\mu/\rho} s_{\rho'/\nu'} = (-1)^{|\nu|} \delta_{\mu, \nu}.
\end{eqnarray}
\end{lemma}

\begin{proof}
We will use the following facts (\cite{Mac}, p. 70).
Let $x=(x_1, x_2, \dots)$ and $y=(y_1, y_2, \dots)$ be two sets of variables
then one has
\begin{eqnarray}
&& \sum_{\mu} s_{\mu}(x)s_{\mu}(y) = \frac{1}{\prod_{i,j} (1 - x_iy_j)}, \\
&& \sum_{\mu} s_{\mu'}(x)s_{\mu}(y) = \prod_{i,j} (1 + x_iy_j),
\end{eqnarray}
and
\begin{eqnarray} \label{eqn:Orth1}
&& \sum_{\mu} s_{\mu/\nu}(x)s_{\mu}(y)
= \sum_{\mu, \rho} c^{\mu}_{\nu\rho}s_{\rho}(x)s_{\mu}(y)
= \sum_{\rho} s_{\rho}(x) s_{\nu}(y)s_{\rho}(y).
\end{eqnarray}
In the same fashion,
we also have
\begin{eqnarray} \label{eqn:Orth2}
&& \sum_{\mu} s_{\mu'/\nu'}(x)s_{\mu}(y)
= \sum_{\mu, \rho} c^{\mu}_{\nu\rho}s_{\rho'}(x)s_{\mu}(y)
= \sum_{\rho} s_{\rho'}(x) s_{\nu}(y)s_{\rho}(y).
\end{eqnarray}
Hence we have
\begin{eqnarray*}
&& \sum_{\mu}\left( \sum_{\rho} (-1)^{\rho} s_{\mu/\rho}(x) s_{\rho'/\nu'}(x)\right) s_{\mu}(y) \\
& = & \sum_{\rho} (-1)^{\rho}  s_{\rho'/\nu'}(x) \sum_{\mu} s_{\mu/\rho}(x) s_{\mu}(y) \\
& = & \sum_{\rho} (-1)^{|\rho|} s_{\rho'/\nu'}(x)\sum_{\theta} s_{\theta}(x)s_{\rho}(y)s_{\theta}(y) \\
& = & \sum_{\theta} s_{\theta}(x) s_{\theta}(y)\sum_{\rho} s_{\rho'/\nu'}(x)s_{\rho}(-y) \\
& = & \sum_{\theta} s_{\theta}(x) s_{\theta}(y)\sum_{\eta} s_{\eta'}(x)s_{\eta}(-y) s_{\nu}(-y) \\
& = &  \frac{1}{\prod_{i,j} (1 - x_iy_j)}\cdot \prod_{i,j} (1 + x_iy_j) \cdot s_{\nu}(-y) \\
& = & s_{\nu}(-y) = (-1)^{|\nu|} s_{\nu}(y).
\end{eqnarray*}
This completes the proof.
\end{proof}

\begin{example}
When $\mu = (n)$ and $\nu = \emptyset$,
then $\rho = (i)$ for $0 \leq i \leq n$.
By (\ref{eqn:SkewDet})
\begin{align*}
s_{\mu/\rho} & = h_{m-i}, & s_{\rho/\nu} & = e_i, \\
\end{align*}
hence (\ref{eqn:Orthogonality}) becomes the well-known identity:
$$\sum_{i=0}^n (-1)^i h_{n-i} e_i = \delta_{n0}.$$
\end{example}

\section{A Reformulation of $\cW_{\mu, \nu}$}
\label{sec:RefW}

\subsection{Some formulas for $\cW_{\mu}$}

Recall
\begin{eqnarray}
&& \cW_{\mu} = q^{\kappa_{\mu}/4}\prod_{1 \leq i < j \leq l(\mu)}
\frac{[\mu_i - \mu_j + j - i]}{[j-i]}
\prod_{i=1}^{l(\mu)} \prod_{v=1}^{\mu_i} \frac{1}{[v-i+l(\mu)]}.
\end{eqnarray}

The following results have been proved in \cite{Zho1}.

\begin{proposition}
We have
\begin{eqnarray}
\cW_{\mu}(q)
& = & q^{-|\mu|/2}s_{\mu}(1, q^{-1}, q^{-2}, \dots) \label{eqn:Wmu1} \\
& = & (-1)^{|\mu|} q^{\kappa_{\mu}/2+|\mu|/2}s_{\mu}(1, q, q^2, \dots). \label{eqn:Wmu2}
\end{eqnarray}
In particular,
\begin{eqnarray} \label{eqn:Wmuqq}
&& q^{-\kappa_{\mu}/2}\cW_{\mu}(q) = \cW_{\mu}(q^{-1}).
\end{eqnarray}
\end{proposition}

\begin{proof}
We have already proved \cite{Zho1}:
\begin{eqnarray} \label{eqn:Wmu}
&& \cW_{\nu}(q) = \frac{q^{\kappa_{\nu}/4}}{\prod_{e \in \nu}
(q^{h(e)/2} - q^{-h(e)/2})}
\end{eqnarray}
We use the following identities from Macdonald's book \cite{Mac}:
$$s_{\mu}(1, q, q^2, \dots) = \frac{q^{n(\mu)}}{\prod_{e \in \mu} ( 1 - q^{h(e)})},$$
and the following identity proved in \cite{Zho1}:
\begin{eqnarray} \label{eqn:hook}
&& \sum_{e\in \mu} h(e) = \frac{1}{2} \kappa_{\mu} + 2 n(\mu) + |\mu|.
\end{eqnarray}
It follows
\begin{eqnarray*}
&& (-1)^{|\mu|} q^{\kappa_{\mu}/2+|\mu|/2}s_{\mu}(1, q, q^2, \dots)
= (-1)^{|\mu|} q^{\kappa_{\mu}/2+|\mu|/2} \cdot
 \frac{q^{n(\mu)}}{\prod_{e \in \mu} ( 1 - q^{h(e)})} \\
& = & \frac{q^{\kappa_{\mu}/2+|\mu|/2+n(\mu)-\frac{1}{2}\sum_{e \in \mu} h(e)}}{\prod_{e \in \mu}
 (q^{h(e)/2} - q^{-h(e)/2})}
= \frac{q^{\kappa_{\mu}/4}}{\prod_{e \in \mu}
 (q^{h(e)/2} - q^{-h(e)/2})}\\
& = & \cW_{\mu}(q)
\end{eqnarray*}
The other identity can be proved in the same fashion.
\end{proof}

\begin{remark}
The identity (\ref{eqn:Wmu1}) holds in the region $|q| > 1$,
while the identity (\ref{eqn:Wmu2}) holds in the region $|q| < 1$.
Note (\ref{eqn:Wmuqq}) holds everywhere as an identity of rational functions
by (\ref{eqn:Wmu}) and (\ref{eqn:hook}).
\end{remark}

\subsection{Definition of $\cW_{\mu, \nu}$}

Recall
\begin{eqnarray}
&& \cW_{\mu, \nu} = q^{|\nu|/2} \cW_{\mu} \cdot s_{\nu}(\cE_{\mu}(q,t)),
\end{eqnarray}
where\begin{eqnarray}
&& \cE_{\mu}(q,t) = \prod_{j=1}^{l(\mu)} \frac{1+q^{\mu_j-j}t}{1+q^{-j}t}
\cdot \left(1 + \sum_{n=1}^{\infty}
\frac{t^n}{\prod_{i=1}^n (q^i-1)}\right).
\end{eqnarray}

By the following identity (cf. e.g. \cite{Mac}, p. 27, Example 4):
\begin{eqnarray*}
&& \prod_{i=1}^{\infty} (1 + q^{i-1}t)
= \sum_{n=0}^{\infty} \frac{q^{n(n-1)/2}t^n}{\prod_{i=1}^n (1 - q^i)}
=  \sum_{n=0}^{\infty} \frac{(q^{-1}t)^n}{\prod_{i=1}^n (q^{-i} - 1)},
\end{eqnarray*}
it follows that
\begin{eqnarray}
&& \cE_{\emptyset}(t) = \prod_{i=1}^{\infty} (1+q^{-i}t),
\end{eqnarray}
and so
\begin{eqnarray}
&& \cE_{\mu}(q,t) = \prod_{j=1}^{l(\mu)} \frac{1+q^{\mu_j-j}t}{1+q^{-j}t}
\cdot \prod_{i=1}^{\infty} (1+q^{-i}t)
=  \prod_{j=1}^{\infty} (1+q^{\mu_j-j}t).
\end{eqnarray}
Therefore,
$$s_{\nu}(\cE_{\mu}(q, t)) = s_{\nu}(q^{\mu_1-1}, q^{\mu_2 -2}, \dots).$$

\subsection{Reformulation of $\cW_{\mu, \nu}$ in terms of skew Schur functions}

We now generalize formula (\ref{eqn:Wmu2}).

\begin{theorem}
We have
\begin{eqnarray} \label{eqn:Keynu}
&& s_{\nu}(\cE_{\mu}(q,t))
= (-1)^{|\nu|} q^{\kappa_{\nu}/2} \sum_{\rho} q^{-|\rho|}
\frac{s_{\mu/\rho}(1, q, q^2, \dots)}{s_{\mu}(1, q, q^2, \dots)}
s_{\nu/\rho}(1, q, q^2, \dots).
\end{eqnarray}
and
\begin{eqnarray} \label{eqn:Key}
&& \cW_{\mu, \nu}(q)
= (-1)^{|\mu|+|\nu|}
q^{\frac{\kappa_{\mu}+\kappa_{\nu}+|\mu|+|\nu|}{2}}
\sum_{\rho} q^{-|\rho|} s_{\mu/\rho}(1, q, \dots)s_{\nu/\rho}(1, q, \dots).
\end{eqnarray}
\end{theorem}

The proof will occupy the rest of this subsection.

Since
\begin{eqnarray*}
\cW_{\mu, \nu}(q) & = & q^{|\nu|/2} \cW_{\mu}(q) s_{\nu}(\cE_{\mu}(q,t)) \\
& = & (-1)^{|\mu|} q^{(\kappa_{\mu}+|\mu|+|\nu|)/2}s_{\mu}(1, q, q^2, \dots)
s_{\nu}(\cE_{\mu}(q,t)),
\end{eqnarray*}
(\ref{eqn:Key}) follows easily from (\ref{eqn:Keynu}).

Let $l = l(\rho)$.
For $n > l(\mu)$,
by (\ref{eqn:SkewDet}) we have
\begin{eqnarray*}
&& s_{\mu/\rho}(1, q, q^2, \cdots)
= \det (h_{\mu_i -\rho_j - i + j})_{1 \leq i, j \leq n},
\end{eqnarray*}
where
$$h_n = \frac{1}{\prod_{i=1}^n (1 - q^i)}.$$
This identity holds in both the region $|q| < 1$
and the region $|q| > 1$.
Here we are working in the latter,
we will change to the former below.
Let $f_0 = 1$ and for $a > 0$, $k \geq 0$,
$$f_{a, k}(q) = (1-q^a)(1-q^{a-1})\cdots (1 - q^{a-(k-1)}),$$
then we have
$$h_n = h_j f_{n, n-j}.$$
It follows that
\begin{eqnarray*}
&& s_{\mu/\rho}(1, q, q^2, \cdots)
= \prod_{i=1}^n h_{\mu_i+n-i} \cdot \det (f_{\mu_i+n-i, \rho_j+n-j})_{1 \leq i, j \leq n}.
\end{eqnarray*}
Since we have
\begin{eqnarray*}
f_{a, k} & = & \sum_{j=0}^k (-1)^j e_j(1, q^{-1}, \dots, q^{-(k-1)}) q^{ja},
\end{eqnarray*}
we can modify the columns as follows.
\begin{eqnarray*}
&& s_{\mu/\rho}(1, q, q^2, \cdots) \\
& = & \prod_{i=1}^n h_{\mu_i+n-i} \cdot \prod_{k=1}^{n-l-1} (-1)^k e_k(1, q^{-1}, \dots, q^{-(k-1)}) \\
&& \cdot \begin{vmatrix}
f_{\mu_1+n-1, \rho_1+n-1} & \cdots & f_{\mu_1+n-1, \rho_l+n-l} & q^{(\mu_1+n-1)(n-l-1)} & \cdots & q^{\mu_1+n-1} & 1 \\
f_{\mu_2+n-2, \rho_1+n-1} & \cdots & f_{\mu_2+n-2, \rho_2+n-l} & q^{(\mu_2+n-2)(n-l-1)} & \cdots & q^{\mu_2+n-2} & 1 \\
\cdot & \cdots & \cdot & \cdot & \cdots & \cdot & \cdot \\
\cdot & \cdots & \cdot & \cdot & \cdots & \cdot & \cdot \\
\cdot & \cdots & \cdot & \cdot & \cdots & \cdot & \cdot \\
f_{\mu_{n-1}+1, \rho_1+n-1} & \cdots & f_{\mu_{n-1}+1, \rho_l+n-l} & q^{(\mu_{n-1}+1)(n-l-1)} & \cdots & q^{\mu_{n-1}+1} & 1 \\
f_{\mu_n, \rho_1+n-1} & \cdots & f_{\mu_n, \rho_l+n-l} &  q^{\mu_n(n-l-1)} & \cdots & q^{\mu_l} & 1
\end{vmatrix}.
\end{eqnarray*}

We introduce the following notations.
Denote by $v^j_{\mu}$ the column vector
$$(q^{(\mu_1+n-1)j}, q^{(\mu_2+n-2)j}, \dots, q^{\mu_nj})^t,$$
and denote by
$$[v_1, \dots, v_n]$$
the matrix formed by the column vectors $v_1, \dots, v_n$.
In these notation we then have
\begin{eqnarray*}
&& s_{\mu/\rho}(1, q, q^2, \cdots) \\
& = & \prod_{i=1}^n h_{\mu_i+n-i} \cdot
\prod_{k=1}^{n-l-1} (-1)^k e_k(1, q^{-1}, \dots, q^{-(k-1)}) \\
&& \cdot \det\left[
\sum_{j_1=0}^{\rho_1+n-1} (-1)^{j_1} e_{j_1}(1, q^{-1}, \dots, q^{-(\rho_1+n-2)}) v^{j_1}_{\mu},
\cdots, \right. \\
&& \left. \sum_{j_l=0}^{\rho_l+n-l} (-1)^{j_l}
e_{j_l}(1, q^{-1}, \dots, q^{-(\rho_l+n-l-1)}) v^{j_l}_{\mu},
v_{\mu}^{n-l-1}, \cdots, v_{\mu}^1, 1\right] \\
& = & \prod_{i=1}^n h_{\mu_i+n-i} \cdot \prod_{k=1}^{n-l-1} (-1)^k e_k(1, q^{-1}, \dots, q^{-(k-1)}) \\
&& \cdot \prod_{i=1}^l \sum_{j_i=n-l}^{\rho_i+n-i}
(-1)^{j_i} e_{j_i}(1, q^{-1}, \dots, q^{-(\rho_i+n-i-1)}) \\
&& \cdot
\det\left[v^{j_1}_{\mu}, \cdots,  v^{j_l}_{\mu}, v_{\mu}^{n-l-1}, \cdots, v_{\mu}^1, 1\right].
\end{eqnarray*}
In particular when $\rho = \emptyset$
we have
\begin{eqnarray*}
s_{\mu}(1, q, q^2, \cdots)
& = & \prod_{i=1}^n h_{\mu_i+n-i} \cdot \prod_{k=1}^{n-1} (-1)^k e_k(1, q^{-1}, \dots, q^{-(k-1)}) \\
&& \cdot \det\left[v^{n-1}_{\mu}, v_{\mu}^{n-2}, \cdots, v_{\mu}^1, 1\right].
\end{eqnarray*}
Therefore,
\begin{eqnarray*}
\frac{s_{\mu/\rho}(1, q, q^2, \dots)}{s_{\mu}(1, q, q^2, \dots)}
& = & \prod_{i=1}^l \sum_{j_i=n-l}^{\rho_i+n-i}
\frac{(-1)^{j_i} e_{j_i}(1, q^{-1}, \dots, q^{-(\rho_i+n-i-1)})}
{(-1)^{n-i} e_{n-i}(1, q^{-1}, \dots, q^{-(n-i-1)})} \\
&& \cdot
\frac{\det\left[v^{j_1}_{\mu}, \cdots,  v^{j_l}_{\mu}, v_{\mu}^{n-l-1}, \cdots, v_{\mu}^1, 1\right]}
{\det\left[v^{n-1}_{\mu}, v_{\mu}^{n-2}, \cdots, v_{\mu}^1, 1\right]}.
\end{eqnarray*}
So far we have been working in the region $|q| < 1$.
We now rewrite the above expressions so that it is easy to take the limit $n \to \infty$
in the region $|q| > 1$.
Take $k_i = j_i - (n-l)$, then we have
\begin{eqnarray*}
&& \frac{(-1)^{j_i} e_{j_i}(1, q^{-1}, \dots, q^{-(\rho_i+n-i-1)})}
{(-1)^{n-i} e_{n-i}(1, q^{-1}, \dots, q^{-(n-i-1)})} \\
& = & \frac{(-1)^{k_i+n-l} e_{k_i+n-l}(1, q^{-1}, \dots, q^{-(\rho_i+n-i-1)})}
{(-1)^{n-i} e_{n-i}(1, q^{-1}, \dots, q^{-(n-i-1)})} \\
& = & (-1)^{k_i-l+i} q^{-\sum_{j=n-i}^{\rho_i+n-i-1} j}
\frac{e_{k_i+n-l}(1, q^{-1}, \dots, q^{-(\rho_i+n-i-1)})}
{1\cdot q^{-1} \dots q^{-(\rho_i+n-i-1)}} \\
& = & (-1)^{k_i-l+i} q^{-\sum_{j=n-i}^{\rho_i+n-i-1} j}
e_{\rho_i+l-i-k_i}(1, q^1, \dots, q^{\rho_i+n-i-1}) \\
& = & (-1)^{k_i-l+i} q^{-\sum_{j=n-i}^{\rho_i+n-i-1} j + (\rho_i+l-i-k_i)(\rho_i+n-i-1)} \\
&& \cdot e_{\rho_i+l-i-k_i}(1, q^{-1}, \dots, q^{-(\rho_i+n-i-1)}).
\end{eqnarray*}
In the above we have used the following easy identities:
\begin{eqnarray*}
&& \frac{e_j(x_1, \dots, x_n)}{x_1 \cdots x_n} = e_{n-j}(x_1^{-1}, \dots, x_n^{-1}), \\
&& e_j(ax_1, \dots, ax_n) = a^j e_j(x_1, \dots, x_n).
\end{eqnarray*}

We deal with the exponent of $q$ as follows.
\begin{eqnarray*}
&& -\sum_{j=n-i}^{\rho_i+n-i-1} j + (\rho_i+l-i-k_i)(\rho_i+n-i-1) \\
& = & - \frac{(\rho_i+n-i-1)(\rho_i+n-i)}{2} + \frac{(n-i-1)(n-i)}{2} \\
&& + (\rho_i+l-i-k_i)(\rho_i+n-i-1) \\
& = & - \frac{\rho_i(\rho_i-1)}{2} + i \rho_i
+ (\rho_i + l - i - k_i)(\rho_i - i - 1)
+ (l - i - k_i)n,
\end{eqnarray*}
and
\begin{eqnarray*}
&& \sum_{i=1}^l (l - i - k_i)n = \left(\frac{l(l-1)}{2} - \sum_{i=1}^l k_i\right)n.
\end{eqnarray*}

Note originally we have $0 \leq k_i \leq \rho_i + l - i$.
However using the convention that
$$e_n = 0$$
for $n < 0$,
one can relieve the constraint $k_i \leq \rho_i + l - i$ in the summation.

Suppose
$$\{k_1+n-l, \dots, k_l + n -l, n- l- 1, \dots, 1, 0\}$$
are
all distinct.
Then there is a partition
$\theta(\vec{k}) = (\theta(\vec{k})_1 \geq \theta(\vec{k})_2 \geq \cdots \geq \theta(\vec{k})_n)$
of length $\leq l$
and a permutation $\sigma_{\vec{k}} \in S_l$
such that
$$\theta(\vec{k})_i+n-i = k_{\sigma_{\vec{k}}(i)} + n - l$$
The degree of this partition is:
\begin{eqnarray*}
|\theta(\vec{k})| & = & \sum_{i=1}^l (k_{\sigma_{\vec{k}}(i)}+l-i)
= \sum_{i=1}^l k_i - \frac{l(l-1)}{2}.
\end{eqnarray*}
Therefore,
\begin{eqnarray*}
&& \frac{\det\left[v^{j_1}_{\mu}, \cdots,  v^{j_l}_{\mu}, v_{\mu}^{n-l-1}, \cdots, v_{\mu}^1, 1\right]}
{\det\left[v^{n-1}_{\mu}, v_{\mu}^{n-2}, \cdots, v_{\mu}^1, 1\right]} \\
& = & \frac{\det\left[v^{k_1+n-l}_{\mu}, \cdots,  v^{k_l+n-l}_{\mu}, v_{\mu}^{n-l-1}, \cdots, v_{\mu}^1, 1\right]}
{\det\left[v^{n-1}_{\mu}, v_{\mu}^{n-2}, \cdots, v_{\mu}^1, 1\right]} \\
& = & (-1)^{\sigma_{\vec{k}}} s_{\theta}(q^{\mu_1+n-1}, \dots, q^{\mu_n}) \\
& = & (-1)^{\sigma_{\vec{k}}} q^{n\left(\sum_{i=1}^l k_i - l(l-1)/2\right)}
s_{\theta(\vec{k})}(q^{\mu_1-1}, \dots, q^{\mu_n-n}).
\end{eqnarray*}
When $$\{k_1+n-l, \dots, k_l + n -l, n- l- 1, \dots, 1, 0\}$$
are not all distinct,
define
$$(-1)^{\sigma_{\vec{k}}} s_{\theta(\vec{k})}(q^{\mu_1-1}, \dots, q^{\mu_n-n}) = 0.$$
Putting things together,
\begin{eqnarray*}
&& \frac{s_{\mu/\rho}(1, q, q^2, \dots)}{s_{\mu}(1, q, q^2, \dots)} \\
& = & \sum_{k_1, \dots, k_l}
\prod_{i=1}^l
(-1)^{k_i-l+i} q^{- \frac{\rho_i(\rho_i-1)}{2} + i \rho_i
+ (\rho_i + l - i - k_i)(\rho_i - i - 1)} \\
&& \cdot
e_{\rho_i+l-i-k_i}(1, q^{-1}, \dots, q^{-(\rho_i+n-i-1)}) \cdot (-1)^{\sigma_{\vec{k}}}
s_{\theta(\vec{k})}(q^{\mu_1-1}, \dots, q^{\mu_n-n}).
\end{eqnarray*}
Hence one can take $n \to \infty$ to get:
\begin{eqnarray*}
\frac{s_{\mu/\rho}(1, q, q^2, \dots)}{s_{\mu}(1, q, q^2, \dots)}
& = & \sum_{k_1, \dots, k_l}
\prod_{i=1}^l \left[
(-1)^{k_i-l+i} q^{- \frac{\rho_i(\rho_i-1)}{2} + i \rho_i
+ (\rho_i + l - i - k_i)(\rho_i - i - 1)} \right. \\
&& \left. \cdot e_{\rho_i+l-i-k_i}(1, q^{-1}, \dots)\right] \cdot (-1)^{\sigma_{\vec{k}}}
s_{\theta(\vec{k})}(\cE_{\mu}(q, t)).
\end{eqnarray*}
Now we apply the identity:
\begin{eqnarray*}
&& e_n(1, q^{-1}, \dots) = (-1)^n q^{-n(n-3)/2} e_n(1, q, \dots),
\end{eqnarray*}
to get
\begin{eqnarray*}
&& \frac{s_{\mu/\rho}(1, q, q^2, \dots)}{s_{\mu}(1, q, q^2, \dots)} \\
& = & \sum_{k_1, \dots, k_l} \prod_{i=1}^l \left[
(-1)^{k_i-l+i} q^{- \frac{\rho_i(\rho_i-1)}{2} + i \rho_i
+ (\rho_i + l - i - k_i)(\rho_i - i - 1)} \right. \\
&& \left. \cdot (-1)^{\rho_i+l-i-k_i} q^{-(\rho_i+l-i-k_i)(\rho_i+l-i-k_i-3)/2}
e_{\rho_i+l-i-k_i}(1, q, \dots) \right] \\
&& \cdot (-1)^{\sigma_{\vec{k}}} s_{\theta(\vec{k})}(\cE_{\mu}(q, t)).
\end{eqnarray*}

The exponent of $q$ is
\begin{eqnarray*}
&& \sum_{i=1}^l [- \rho_i(\rho_i-1)/2 + i \rho_i
+ (\rho_i + l - i - k_i)(\rho_i - i - 1)  \\
&& -(\rho_i+l-i-k_i)(\rho_i+l-i-k_i-3)/2] \\
& = & \sum_{i=1}^l (l-i+2\rho_i-k_i+i^2- l^2+2lk_i- k_i^2)/2 \\
& = & \frac{1}{2} \left(2|\rho| + (2l-1)\sum_{i=1}^l k_i - \sum_{i=1}^l k_i^2
- \frac{l(l-1)(2l-1)}{3}\right).
\end{eqnarray*}
Note this is invariant under the permutation of $k_1, \dots, k_l$,
hence it depends only on the partition $\theta(\vec{k})$,
and so one can take
$$k_i = \theta(\vec{k})_i + l - i,$$
then the exponent of $q$ becomes one half of:
\begin{eqnarray*}
&& 2|\rho| + (2l-1) \sum_{i=1}^l (\theta(\vec{k})_i+l-i)
- \sum_{i=1}^l (\theta(\vec{k})_i+l-i)^2 - \frac{1}{3}l(l-1)(2l-1) \\
& = & 2|\rho| + (2\rho - 1) |\theta(\vec{k})| + (2l-1) \cdot \frac{l(l-1)}{2}
- \sum_{i=1}^l (\theta(\vec{k})_i^2 - 2i\theta(\vec{k})_i)
- 2l |\theta(\vec{k})| \\
&& - \sum_{i=1}^l (l-i)^2 - \frac{1}{3}l(l-1)(2l-1) \\
& = & 2 |\rho| - \kappa_{\theta(\vec{k})}.
\end{eqnarray*}

To summarize,
it is possible to rewrite the summation over $k_1, \dots, k_n$ as a summation over all
partitions $\theta$ of length $\leq l$:
\begin{eqnarray*}
&& \frac{s_{\mu/\rho}(1, q, q^2, \dots)}{s_{\mu}(1, q, q^2, \dots)} \\
& = & (-1)^{|\rho|} \sum_{l(\theta) \leq l} q^{|\rho| - \kappa_{\theta}/2}
\sum_{\sigma \in S_l} (-1)^{\sigma} \prod_{i=1}^l e_{\rho_i-i-\theta_{\sigma(i)} + \sigma(i)}(1, q, \dots)
\cdot s_{\theta}(\cE_{\mu}(q, t)) \\
& = & (-1)^{|\rho|} \sum_{l(\theta) \leq l}
q^{|\rho| - \kappa_{\theta}/2}
\det (e_{\rho_i-i-\theta_j + j}(1, q, \dots))_{1 \leq i, j \leq l}
\cdot s_{\theta}(\cE_{\mu}(q, t)) \\
& = & (-1)^{|\rho|} \sum_{l(\theta) \leq l}
q^{|\rho| - \kappa_{\theta}/2}
s_{\rho'/\theta'}(1, q, \dots)
\cdot s_{\theta}(\cE_{\mu}(q, t))  \\
& = & (-1)^{|\rho|} \sum_{\theta}
q^{|\rho| - \kappa_{\theta}/2}
s_{\rho'/\theta'}(1, q, \dots)
\cdot s_{\theta}(\cE_{\mu}(q, t)).
\end{eqnarray*}
Finally we have
\begin{eqnarray*}
&& (-1)^{|\nu|} q^{\kappa_{\nu}/2} \sum_{\rho} q^{-|\rho|}
\frac{s_{\mu/\rho}(1, q, q^2, \dots)}{s_{\mu}(1, q, q^2, \dots)}
s_{\nu/\rho}(1, q, q^2, \dots) \\
& = & (-1)^{|\nu|} q^{\kappa_{\nu}/2} \sum_{\rho} q^{-|\rho|}
(-1)^{|\rho|} \sum_{\theta}
q^{|\rho| - \kappa_{\theta}/2}
s_{\rho'/\theta'}(1, q, \dots) s_{\theta}(\cE_{\mu}(q, t))
s_{\nu/\rho}(1, q,  \dots) \\
& = & (-1)^{|\nu|} q^{\kappa_{\nu}/2}\sum_{\theta}
q^{-\kappa_{\theta}/2} s_{\theta}(\cE_{\mu}(q, t))
\sum_{\rho} (-1)^{|\rho|} s_{\rho'/\theta'}(1, q, \dots)
s_{\nu/\rho}(1, q, q^2, \dots) \\
& = & (-1)^{|\nu|} q^{\kappa_{\nu}/2}\sum_{\theta}
q^{-\kappa_{\theta}/2} s_{\theta}(\cE_{\mu}(q, t))
(-1)^{|\nu|}\delta_{\theta, \nu} \\
& = & s_{\nu}(\cE_{\mu}(q, t)).
\end{eqnarray*}
In the second to last equality we have used (\ref{eqn:Orthogonality}).
This finishes the proof of (\ref{eqn:Keynu}).

\section{Initial Values at $\tau = -1$}
\label{sec:Initial}

Now we study the initial values at $\tau = -1$ of the two sides of (\ref{eqn:Conjecture}).

\subsection{The left-hand side}

When $l(\mu^+) + l(\mu^-) > 2$,
$$G_{\mu^+, \mu^-}(\lambda; -1) = 0;$$
when $l(\mu^+) = 1$ and $l(\mu^-) = 0$,
\begin{eqnarray*}
&& G_{\mu^+, \mu^-}(\lambda; -1) \\
& = &  - \sqrt{-1} \lambda^{-1}  \sum_{g \geq 0} \lambda^{2g}
\int_{\Mbar_{g, 1}} \frac{\lambda_g} {\frac{1}{\mu_1^+}
\left(\frac{1}{\mu^+_1} - \psi_1\right)}
\frac{\prod_{a=1}^{\mu^+_1-1} \left(-\mu_1^++ a \right)}{\mu_1^+ \cdot \mu_1^+!} \\
& = &  (-1)^{\mu^+_1} \sqrt{-1} \cdot \frac{1}{2\mu_1^+\sin (\mu_1^+\lambda/2)}
= \frac{(-1)^{\mu_1^+-1}}{q^{\mu_1^+/2} - q^{-\mu_1^+/2}} \cdot \frac{p_{\mu_1^+}}{\mu_1^+};
\end{eqnarray*}
the case of $l(\mu^+) = 0$ and $l(\mu^-) = 1$ is similar;
when $l(\mu^+) = l(\mu^-) = 1$,
\begin{eqnarray*}
G_{\mu^+, \mu^-}(\lambda; -1)
& = &  \lim_{\tau \to -1}  \sum_{g \geq 0} \lambda^{2g} \int_{\Mbar_{g, 2}}
\frac{\Lambda_{g}^{\vee}(1)\Lambda^{\vee}_{g}(\tau)\Lambda_{g}^{\vee}(-1 - \tau)}
{\frac{1}{\mu_1^+} \left(\frac{1}{\mu^+_1} - \psi_1\right)
\cdot \frac{\tau}{\mu^-_1} \left(\frac{\tau}{\mu^-_1} - \psi_{2}\right)} \\
&& \cdot \tau(1+\tau) \cdot
\frac{\prod_{a=1}^{\mu^+_1-1} \left(\mu_1^+ \tau + a \right)}{\mu_1^+ \cdot \mu_1^+!}
\cdot \frac{\prod_{a=1}^{\mu^-_1-1} \left(\frac{\mu_1^-}{\tau} + a \right)}{\mu_1^- \cdot \mu_1^-!}.
\end{eqnarray*}
One needs to consider the $g=0$ term and the $g>0$ terms separately.
In the second case, the limit is zero while in first case,
by our convention:
$$\int_{\Mbar_{0, 2}}
\frac{\Lambda_{0}^{\vee}(1)\Lambda^{\vee}_{0}(\tau)\Lambda_{0}^{\vee}(-1 - \tau)}
{\frac{1}{\mu_1^+} \left(\frac{1}{\mu^+_1} - \psi_1\right)
\cdot \frac{\tau}{\mu^-_1} \left(\frac{\tau}{\mu^-_1} - \psi_{2}\right)}
= \frac{(\mu_1^+)^2(\frac{\mu_1^-}{\tau})^2}{\mu_1^+ + \frac{\mu_1^-}{\tau}},$$
hence when $\mu_1^+ \neq \mu_1^-$,
the limit is zero,
when $\mu_1^+ = \mu_1^-$,
the limit is:
\begin{eqnarray*}
&& \lim_{\tau \to -1}
\frac{(\mu_1^+)^2(\frac{\mu_1^-}{\tau})^2}{\mu_1^+ + \frac{\mu_1^-}{\tau}}
\cdot
\tau(1+\tau) \cdot
\frac{\prod_{a=1}^{\mu^+_1-1} \left(\mu_1^+ \tau + a \right)}{\mu_1^+ \cdot \mu_1^+!}
\cdot \frac{\prod_{a=1}^{\mu^-_1-1} \left(\frac{\mu_1^-}{\tau} + a \right)}{\mu_1^-\cdot\mu_1^-!}
= \frac{1}{\mu_1^+}.
\end{eqnarray*}
To summarize,
the initial value is:
\begin{eqnarray*}
&& G(\lambda;p(x^+),p(x^-);-1)^{\bullet} \\
& = & \exp \left(\sum_{n \geq 1} \frac{(-1)^{n-1}}{q^{n/2} - q^{-n/2}}\frac{p_n(x^+)}{n}
+ \sum_{n \geq 1} \frac{(-1)^{n-1}}{q^{n/2} - q^{-n/2}}\frac{p_n(x^-)}{n}
+ \sum_{n \geq 1} \frac{p_n(x^+)p_n(x^-)}{n}\right) \\
& = & \prod_{i,j=1}^{\infty} \frac{1}{(1 + q^{i-1/2}x^+_j)(1+q^{i-1/2}x^-_j)}
\prod_{j,k} \frac{1}{1-x_j^+x^-_k} \\
& = & \sum_{\rho^+} s_{\rho^+}(-q^{1/2}, -q^{3/2}, \dots)s_{\rho^+}(x^+)
\cdot \sum_{\rho} s_{\rho}(x^+)s_{\rho}(x^-) \\
&& \cdot \sum_{\nu^-} s_{\rho^-}(-q^{1/2}, -q^{3/2}, \dots)s_{\rho^-}(x^-) \\
& = & \sum_{\nu^{\pm}, \rho, \rho^{\pm}}
 s_{\rho^+}(-q^{1/2}, -q^{3/2}, \dots)c_{\rho^+\rho}^{\nu^+}s_{\nu^+}(x^+)
\cdot c_{\rho^-\rho}^{\nu^-} s_{\nu^-}(x^-) s_{\rho^-}(-q^{1/2}, -q^{3/2}, \dots) \\
& = & \sum_{\rho, \nu^{\pm}}
s_{\nu^+/\rho}(-q^{1/2}, -q^{3/2}, \dots)s_{\nu^-/\rho}(-q^{1/2}, -q^{3/2}, \dots)
\cdot s_{\nu^+}(x^+) s_{\nu^-}(x^-).
\end{eqnarray*}

\subsection{The right-hand side}

\begin{eqnarray*}
&& R(-1; p(x^+), p(x^-)) \\
& = & \sum_{|\mu^{\pm}|=|\nu^{\pm}|}
\frac{\chi_{\nu^+}(\mu^+)}{z_{\mu^+}}
\frac{\chi_{\nu^-}(\mu^-)}{z_{\mu^-}}
e^{-\sqrt{-1}(\kappa_{\nu^+}+ \kappa_{\nu^-})\lambda/2} \cW_{\nu^+, \nu^-}(q)
p^+_{\mu^+}p^-_{\mu^-} \\
& = & \sum_{\nu^{\pm}} s_{\nu^+}(x^+) q^{-\kappa_{\nu^+}/2}
\cW_{\nu^+, \nu^-}(q) q^{-\kappa_{\nu^-}/2}s_{\nu^-}(x^-) \\
& = & \sum_{\nu^{\pm}} s_{\nu^+}(x^+)
s_{\nu^-}(x^-)
(-1)^{|\nu^+|+|\nu^-|}q^{(|\nu^+|+|\nu^-|)/2}
 \sum_{\rho} q^{-|\rho|} s_{\nu^+/\rho}(1, q, \dots)s_{\nu^-/\rho}(1, q, \dots) \\
& = & \sum_{\nu^{\pm}} s_{\nu^+}(x^+)
s_{\nu^-}(x^-)
 \sum_{\rho} s_{\nu^+/\rho}(-q^{1/2}, -q^{3/2}, \dots)s_{\nu^-/\rho}(-q^{1/2}, -q^{3/2}, \dots).
\end{eqnarray*}

So we have proved:

\begin{theorem}
$$R(-1; p(x^+), p(x^-)) = G(-1; p(x^+), p(x^-))^{\bullet}.$$
\end{theorem}

\section{Initial Values at $\tau = 1$}

We study in this section the initial values at $\tau = 1$ of the two sides of (\ref{eqn:Conjecture}).
We need the following
\begin{lemma}
Let $p_{\mu} = \prod_i (p_i(x^+) + p_i(x^-))^{m_i(\mu)}.$
Then we have
\begin{eqnarray*}
&& \sum_{\mu^+\cup \mu^- = \mu} \frac{z_{\mu}}{z_{\mu^+} \cdot z_{\mu^-}}
p_{\mu^+}(x^+) p_{\mu^-}(x^-)
= p_{\mu}.
\end{eqnarray*}
\end{lemma}

\begin{proof}
\begin{eqnarray*}
&& \sum_{\mu^+\cup \mu^- = \mu} \frac{z_{\mu}}{z_{\mu^+} \cdot z_{\mu^-}}
p_{\mu^+}(x^+) p_{\mu^-}(x^-) \\
& = & \prod_i \sum_{j=0}^{m_i(\mu)} \frac{m_i(\mu)!}{j!(m_i(\mu)-j)!} p_i(x^+)^jp_i(x^-)^{m_i(\mu)-j} \\
& = & \prod_i (p_i(x^+) + p_i(x^-))^{m_i(\mu)} \\
& = & p_{\mu}.
\end{eqnarray*}
\end{proof}

Now notice that
\begin{eqnarray*}
&& G_{\mu^+, \mu^-}(\lambda; 1) \\
& = &  -\frac{\sqrt{-1}^{l(\mu^+)+l(\mu^-)}}{z_{\mu^+} \cdot z_{\mu^-}}
\lambda^{l(\mu^+)+l(\mu^-)-2} \\
&& \cdot \sum_{g \geq 0} \lambda^{2g} \int_{\Mbar_{g, l(\mu^+)+l(\mu^-)}}
\frac{\Lambda_{g}^{\vee}(1)\Lambda^{\vee}_{g}(1)\Lambda_{g}^{\vee}(-2)}
{\prod_{i=1}^{l(\mu^+)} \frac{1}{\mu_i^+}
\left(\frac{1}{\mu^+_i} - \psi_i\right)
\prod_{j=1}^{l(\mu^-)} \frac{1}{\mu^-_j} \left(\frac{1}{\mu^-_j} - \psi_{j+l(\mu^+)}\right)} \\
&& \cdot \left[1(1+1)\right]^{l(\mu^+)+l(\mu^-)-1} \cdot
\prod_{i=1}^{l(\mu^+)} \frac{\prod_{a=1}^{\mu^+_i-1}(\mu_i^+  + a)}{\mu_i^+!}
\cdot \prod_{i=1}^{l(\mu^-)} \frac{\prod_{a=1}^{\mu^-_i-1} (\mu_i^- + a)}{\mu_i^-!} \\
& = & \frac{z_{\mu^+\cup \mu^-}}{z_{\mu^+} \cdot z_{\mu^-}}
G_{\mu^+ \cup \mu^-}(\lambda;1).
\end{eqnarray*}
Hence
\begin{eqnarray*}
G(\lambda;1;p^+, p^-)
& = & \sum_{(\mu^+, \mu^-)\in \cP^2_+} G_{\mu^+ \cup \mu^-}(\lambda;1)
\frac{z_{\mu^+\cup \mu^-}}{z_{\mu^+} \cdot z_{\mu^-}} p_{\mu^+}^+p_{\mu^-}^-
= \sum_{|\mu|>0} G_{\mu}(\lambda;1)p_{\mu}, \\
G^{\bullet}(\lambda;1;p^+, p^-)
& = & \exp\left(\sum_{|\mu|>0} G_{\mu}(\lambda;1)p_{\mu}\right)
= \sum_{|\mu|=|\nu|\geq 0} \frac{\chi_{\nu}(\mu)}{z_{\mu}}\cW_{\nu} q^{\kappa_{\nu}/2} p_{\mu} \\
& = & \sum_{|\mu|=|\nu|\geq 0} \sum_{\mu^+\cup\mu^-=\mu}
\frac{\chi_{\nu}(\mu)}{z_{\mu^+}\cdot z_{\mu^-}}
\cW_{\nu} q^{\kappa_{\nu}/2} p_{\mu^+}^+p_{\mu^-}^-.
\end{eqnarray*}

On the other hand, the right-hand side is given by:
\begin{eqnarray*}
&& R(1; p(x^+), p(x^-)) \\
& = & \sum_{|\mu^{\pm}|=|\nu^{\pm}|}
\frac{\chi_{\nu^+}(\mu^+)}{z_{\mu^+}}
\frac{\chi_{\nu^-}(\mu^-)}{z_{\mu^-}}
q^{(\kappa_{\nu^+}+ \kappa_{\nu^-})/2} \cW_{\nu^+, \nu^-}
p^+_{\mu^+}p^-_{\mu^-}.
\end{eqnarray*}
Therefore,
by comparing the coefficients of $p^+_{\mu^+}p^-_{\mu^-}$,
we get:
\begin{eqnarray*}
&&  \sum_{|\nu|\geq 0} \sum_{\mu=\mu^+\cup\mu^-}
\frac{\chi_{\nu}(\mu)}{z_{\mu^+}\cdot z_{\mu^-}}
\cW_{\nu} q^{\kappa_{\nu}/2} \\
& = & \sum_{|\nu^{\pm}| \geq 0}
\frac{\chi_{\nu^+}(\mu^+)}{z_{\mu^+}}
\frac{\chi_{\nu^-}(\mu^-)}{z_{\mu^-}}
q^{(\kappa_{\nu^+}+ \kappa_{\nu^-})/2} \cW_{\nu^+, \nu^-}.
\end{eqnarray*}
By the orthogonality relations:
\begin{eqnarray*}
&& \sum_{\mu} \frac{\chi_{\nu}(\mu)\chi_{\eta}(\mu)}{z_{\mu}} = \delta_{\nu, \eta},
\end{eqnarray*}
we then easily get

\begin{theorem}
Assuming (\ref{eqn:Conjecture}) we have
\begin{eqnarray}
&\cW_{\nu^+, \nu^-}
= q^{-(\kappa_{\nu^+}+ \kappa_{\nu^-})/2} \sum_{\nu,\mu^+,\mu^-}
\frac{\chi_{\nu}(\mu^+\cup \mu^-)\chi_{\nu^+}(\mu^+)\chi_{\nu^-}(\mu^-)}{z_{\mu^+}\cdot z_{\mu^-}}
\cW_{\nu} q^{\kappa_{\nu}/2}.
\end{eqnarray}
\end{theorem}

\section{Examples}
\label{sec:Examples}

\subsection{Examples of $\cW_{\mu, \nu}$}

Using (\ref{eqn:Key}) it is straightforward to find:
\begin{eqnarray*}
&& \cW_{(1), (1)} - \cW_{(1)}\cW_{(1)} = 1, \\
&& \cW_{(2), (1)} - \cW_{(2)}\cW_{(1)} = \frac{q}{q^{1/2} - q^{-1/2}}, \\
&& \cW_{(11), (1)} - \cW_{(11)}\cW_{(1)}
= \frac{q^{-1}}{q^{1/2} - q^{-1/2}}, \\
&& \cW_{(3), (1)} - \cW_{(3)}\cW_{(1)}
= \frac{q^{5/2}}{(q^{1/2} - q^{-1/2})(q-q^{-1})}, \\
&& \cW_{(21), (1)} - \cW_{(2)}\cW_{(1)}
= \frac{1}{(q^{1/2} - q^{-1/2})^2}, \\
&& \cW_{(1^3), (1)} - \cW_{(1^3)}\cW_{(1)}
= \frac{q^{5/2}}{(q^{1/2} - q^{-1/2})(q-q^{-1})}, \\
&& \cW_{(2), (11)} - \cW_{(2)}\cW_{(11)}
= \frac{1}{(q^{1/2} - q^{-1/2})^2}, \\
&& \cW_{(11), (11)} - \cW_{(11)}\cW_{(11)}
= \frac{q^{-2}(q^{3/2} - q^{-3/2})}{(q^{1/2} - q^{-1/2})(q-q^{-1})}, \\
&& \cW_{(3), (11)} - \cW_{(3)}\cW_{(11)} = \frac{q^{3/2}}{(q^{1/2} - q^{-1/2})^2(q-q^{-1})}, \\
&& \cW_{(21), (11)} - \cW_{(21)}\cW_{(11)}
= \frac{q^{-1}(q^{3/2} + q^{-3/2})}{(q^{1/2} - q^{-1/2})^2(q-q^{-1})}, \\
&& \cW_{(1^3), (11)} - \cW_{(1^3)}\cW_{(11)}
= \frac{q^{-5/2}}{(q^{1/2} - q^{-1/2})^2(q-q^{-1})}
+ \frac{q^{-5}}{(q^{1/2} - q^{-1/2})(q-q^{-1})}, \\
&& \cW_{(2), (2)} - \cW_{(2)}\cW_{(2)}
= \frac{q^2(q^{3/2} + q^{-3/2})}{(q^{1/2} - q^{-1/2})(q-q^{-1})}, \\
&& \cW_{(3), (2)} - \cW_{(3)}\cW_{(2)}
= \frac{q^5}{(q^{1/2} - q^{-1/2})(q-q^{-1})}
+ \frac{q^{5/2}}{(q^{1/2} - q^{-1/2})^2(q-q^{-1})}, \\
&& \cW_{(1^3), (2)} - \cW_{(1^3)}\cW_{(2)}
= \frac{q^{-3/2}}{(q^{1/2} - q^{-1/2})^2(q-q^{-1})}.
\end{eqnarray*}

\subsection{The prediction for $G_{(n), (1)}$}

By (\ref{eqn:Conjecture}) we have
\begin{eqnarray} \label{eqn:Predictionn1}
&& G_{(n), (1)}(\lambda;\tau)
= \sum_{|\nu^+|=n} \frac{\chi_{\nu^+}(C_{(n)})}{n}
e^{\sqrt{-1}\tau \kappa_{\nu^+}\lambda/2}
(\cW_{\nu^+, (1)} - \cW_{\nu^+}\cW_{(1)}).
\end{eqnarray}
For example,
when $n =1$,
\begin{eqnarray*}
G_{(1), (1)}(\tau) = \cW_{(1, 1)} - \cW_{(1)}^2
= \frac{q+q^{-1}-1}{(q^{1/2}- q^{-1/2})^2} - \frac{1}{(q^{1/2} - q^{-1/2})^2}
= 1,
\end{eqnarray*}
or equivalently,
\begin{eqnarray*}
&& \sum_{g \geq 0} \lambda^{2g} \int_{\Mbar_{g, 2}}
\frac{\Lambda_{g}^{\vee}(1)\Lambda^{\vee}_{g}(\tau)\Lambda_{g}^{\vee}(-1 - \tau)}
{(1 - \psi_1)(\tau - \psi_2)}
= \frac{1}{\tau(\tau+1)}.
\end{eqnarray*}
When $n=2$,
\begin{eqnarray*}
&& G_{(2), (1)}(\tau) \\
& = &
\frac{\cW_{(2), (1)} - \cW_{(2)}\cW_{(1)}}{2}e^{\sqrt{-1}\tau\lambda}
- \frac{\cW_{(1,1), (1)} - \cW_{(1,1)}\cW_{(1)}}{2}e^{-\sqrt{-1}\tau\lambda} \\
& = & \frac{q}{q^{1/2} - q^{-1/2}} \frac{e^{\sqrt{-1}\tau\lambda}}{2}
- \frac{q^{-1}}{q^{1/2} - q^{-1/2}} \frac{e^{-\sqrt{-1}\tau\lambda}}{2}\\
& = & \frac{\sin [(\tau + 1)\lambda]}{2\sin(\lambda/2)}.
\end{eqnarray*}
When $n =3$,
\begin{eqnarray*}
G_{(3), (1)}(\tau)
& = & -\frac{\cos [(3\tau + 5/2)\lambda]}{6\sin(\lambda/2)\sin\lambda}
+ \frac{1}{12\sin^2(\lambda/2)} \\
& = & \frac{\sin[(3\tau+2)\lambda/2]\sin[(3\tau+3)\lambda/2]}{3\sin(\lambda/2)\sin\lambda}.
\end{eqnarray*}
These result provide evidence for (\ref{eqn:Cn1}) below.

Recall the following conjecture made in \cite{Zho1}:
\begin{eqnarray} \label{eqn:Conjecture-}
&& C_{n,1}(\lambda;p)
= n^2 \cdot n!  \prod_{i=0}^{n}
\frac{\sin [(np+i)\lambda/2]}{(np+i)\sin [(i+1)\lambda/2]},
\end{eqnarray}

\begin{proposition}
Assuming (\ref{eqn:Conjecture-}),
we have
\begin{eqnarray} \label{eqn:Cn1}
&& G_{(n), (1)}(x, y) = \frac{1}{n} \prod_{i=1}^{n-1}
\frac{\sin [(n\tau+i+1)\lambda/2]}{\sin (i\lambda/2)}.
\end{eqnarray}
\end{proposition}

\begin{proof}
We have
\begin{eqnarray*}
&& G_{(n), (1)}(x, y) \\
& = & \frac{1}{n} \sum_{g \geq 0} \lambda^{2g} \int_{\Mbar_{g, 2}}
\frac{\Lambda_{g}^{\vee}(x)\Lambda^{\vee}_{g}(y)\Lambda_{g}^{\vee}(-x - y)}
{\frac{x}{n}\left(\frac{x}{n} - \psi_1\right) y(y - \psi_2)}
\cdot xy(x+y) \cdot \frac{\prod_{a=1}^{n-1} (ax + ny)}{n!x^{n-1}} \\
& = & \frac{\tau(\tau+1)\prod_{a=1}^{n-1} (a + n\tau)}{n \cdot n!} \cdot
\sum_{g \geq 0} \lambda^{2g} \int_{\Mbar_{g, 2}}
\frac{\Lambda_{g}^{\vee}(1)\Lambda^{\vee}_{g}(\tau)\Lambda_{g}^{\vee}(-1 - \tau)}
{\frac{1}{n}\left(\frac{1}{n} - \psi_1\right)\tau(\tau - \psi_2)} \\
& = & \frac{\tau(\tau+1)\prod_{a=1}^{n-1} (a + n\tau)}{n \cdot n!\tau^3} \cdot
\sum_{g \geq 0} \lambda^{2g} \int_{\Mbar_{g, 2}}
\frac{\Lambda_{g}^{\vee}(\frac{1}{\tau})\Lambda^{\vee}_{g}(1)\Lambda_{g}^{\vee}(-1 - \frac{1}{\tau})}
{\frac{1}{n\tau}\left(\frac{1}{n\tau} - \psi_1\right)(1 - \psi_2)}
\end{eqnarray*}
where in the last two equalities we have applied Lemma \ref{lm:Degree}.
By (\ref{eqn:Conjecture-}) one gets:
\begin{eqnarray*}
&& G_{(n), (1)}(x, y) \\
& = & \frac{\tau(\tau+1)\prod_{a=1}^{n-1} (a + n\tau)}{n \cdot n!\tau^3}
\cdot (n\tau)^2 (n\tau)! \prod_{i=0}^{n\tau}
\frac{\sin [(n\tau \cdot \frac{1}{\tau}+i)\lambda/2]}{(n\tau \cdot \frac{1}{\tau}+i)\sin [(i+1)\lambda/2]} \\
& = & \frac{1}{n} \prod_{i=1}^{n-1}
\frac{\sin [(n\tau+i+1)\lambda/2]}{\sin (i\lambda/2)}.
\end{eqnarray*}
\end{proof}

\subsection{The prediction for $G_{(n),(2)}$}

Similarly,
from (\ref{eqn:Conjecture}) one can get:
\begin{equation} \label{eqn:Predictionn2}
\begin{split}
& G_{(n), (2)}(\lambda;\tau) \\
= & \frac{1}{2} \sum_{|\nu^+|=n} \frac{\chi_{\nu^+}(C_{(n)})}{n}
e^{\sqrt{-1}(\tau \kappa_{\nu^+} + \frac{2}{\tau})\lambda/2}
(\cW_{\nu^+, (2)} - \cW_{\nu^+}\cW_{(2)}) \\
- & \frac{1}{2} \sum_{|\nu^+|=n} \frac{\chi_{\nu^+}(C_{(n)})}{n}
e^{\sqrt{-1}(\tau \kappa_{\nu^+} - \frac{2}{\tau})\lambda/2}
(\cW_{\nu^+, (1,1)} - \cW_{\nu^+}\cW_{(1,1)}).
\end{split} \end{equation}
The $n=1$ case is equivalent to the prediction for $G_{(2), (1)}(\tau)$.
In the $n=2$ case one gets
\begin{eqnarray*}
G_{(2), (2)}(\lambda;\tau)
& = & -\frac{\cos[(\tau + \frac{1}{\tau} + 2)\lambda]\cos(3\lambda/2)}{4\sin(\lambda/2)\sin\lambda}
+ \frac{\cos [(\tau-\frac{1}{\tau})\lambda]}{8\sin^2(\lambda/2)} \\
& = & \frac{\sin[2(\tau + \frac{1}{2})(\frac{1}{\tau}+\frac{1}{2})\lambda]}{4\sin(\lambda/2)}
+ \frac{\sin [(\tau + \frac{1}{2})\lambda]\sin[(\frac{1}{\tau}+\frac{1}{2})\lambda]}
{4\sin^2(\lambda/2)}.
\end{eqnarray*}
In the $n=3$ one gets
\begin{eqnarray*}
&& G_{(3), (2)}(\tau) \\
& = & - \frac{\cos [(3\tau+\frac{1}{\tau}+5)\lambda]}{12\sin(\lambda/2)\sin\lambda}
- \frac{\sin[(3\tau + \frac{1}{\tau}+\frac{5}{2})\lambda]}
{24\sin^2(\lambda/2)\sin\lambda} \\
&& + \frac{\sin[(3\tau - \frac{1}{\tau}+\frac{3}{2})\lambda]}
{24\sin^2(\lambda/2)\sin\lambda}
+ \frac{\sin[(\frac{1}{\tau}+1)\lambda]\cos(3\lambda/2)}
{12\sin^2(\lambda/2)\sin\lambda} \\
& = & \frac{1}{24\sin^2(\lambda/2)\sin\lambda}
\left(-\sin[(3\tau+\frac{1}{\tau}+\frac{11}{2})\lambda]
+ \sin [(3\tau+\frac{1}{\tau}+\frac{9}{2})\lambda] \right. \\
&& \left. - \sin[(3\tau+\frac{1}{\tau}+\frac{5}{2})\lambda]
+ \sin [(3\tau-\frac{1}{\tau}+\frac{3}{2})\lambda]
+ \sin[(\frac{1}{\tau}+\frac{5}{2})\lambda]
+ \sin [(\frac{1}{\tau}-\frac{1}{2})\lambda]\right).
\end{eqnarray*}

\subsection{The prediction for $G_{(n),(1,1)}$}
From (\ref{eqn:Conjecture}) one can get:
\begin{equation} \label{eqn:Predictionn11}
\begin{split}
& G_{(n), (1,1)}(\lambda;\tau) \\
= & \frac{1}{2} \sum_{|\nu^+|=n} \frac{\chi_{\nu^+}(C_{(n)})}{n}
e^{\sqrt{-1}(\tau \kappa_{\nu^+} + \frac{2}{\tau})\lambda/2}
(\cW_{\nu^+, (2)} - \cW_{\nu^+}\cW_{(2)}) \\
+ & \frac{1}{2} \sum_{|\nu^+|=n} \frac{\chi_{\nu^+}(C_{(n)})}{n}
e^{\sqrt{-1}(\tau \kappa_{\nu^+} - \frac{2}{\tau})\lambda/2}
(\cW_{\nu^+, (1,1)} - \cW_{\nu^+}\cW_{(1,1)}) \\
- & \sum_{|\nu^+|=n} \frac{\chi_{\nu^+}(C_{(n)})}{n}
e^{\sqrt{-1}\tau \kappa_{\nu^+}\lambda/2}
(\cW_{\nu^+, (1)} - \cW_{\nu^+}\cW_{(1)}) \cW_{(1)}.
\end{split} \end{equation}

When $n=1$ one gets
\begin{eqnarray*}
G_{(1), (1,1)}(\tau)
& = & \frac{\cos [(\frac{1}{\tau} + 1)\lambda]-1}{2\sqrt{-1}\sin(\lambda/2)}
= \frac{\sqrt{-1} \sin^2[(\frac{1}{\tau}+1)\lambda/2]}{\sin(\lambda/2)}.
\end{eqnarray*}
On the other hand,
we have
\begin{eqnarray*}
C_{(1,1), (1)}(\tau)
& = & \frac{\cos [(\tau + 1)\lambda]-1}{2\sqrt{-1}\sin(\lambda/2)}
= \frac{\sqrt{-1} \sin^2[(\tau+1)\lambda/2]}{\sin(\lambda/2)}.
\end{eqnarray*}
This is compatible with (\ref{eqn:Inverse}).
When $n=2$,
\begin{eqnarray*}
&& G_{(2), (1,1)}(\lambda;\tau) \\
& = & \frac{\sin [(\tau + \frac{1}{\tau} + 2)\lambda]\cos(3\lambda/2)}
{4\sqrt{-1}\sin(\lambda/2)\sin\lambda}
+ \frac{\sin[(\tau - \frac{1}{\tau})\lambda]}{8\sqrt{-1}\sin^2(\lambda/2)}
- \frac{\sin [(\tau+1)\lambda]}{4\sqrt{-1}\sin^2(\lambda/2)}\\
& = & \frac{\sqrt{-1}\sin [(\tau+1)\lambda]\sin^2[(\frac{1}{\tau} + 1)\lambda/2]}
{2\sin^2(\lambda/2)}.
\end{eqnarray*}
When $n=3$
\begin{eqnarray*}
& & G_{(3), (1,1)}(\lambda;\tau) \\
& = & -\frac{\sqrt{-1}\sin[(3\tau+\frac{1}{\tau}+5)\lambda]}
{12\sin(\lambda/2)\sin\lambda}
+ \frac{\sqrt{-1}\cos[(3\tau+\frac{1}{\tau}+\frac{5}{2})\lambda]}
{24\sin^2(\lambda/2)\sin\lambda}
+ \frac{\sqrt{-1}\cos[(3\tau-\frac{1}{\tau}+\frac{3}{2})\lambda]}
{24\sin^2(\lambda/2)\sin\lambda} \\
&& -  \frac{\sqrt{-1}\cos[(\frac{1}{\tau}+1)\lambda]\cos(3\lambda/2)}
{12\sin^2(\lambda/2)\sin\lambda}
- \frac{\sqrt{-1}\cos[(3\tau+\frac{5}{2})\lambda]}
{12\sin^2(\lambda/2)\sin\lambda}
+ \frac{\sqrt{-1}}{12\sin^3(\lambda/2)}.
\end{eqnarray*}

\end{document}